\documentstyle[12pt]{article}

  \newcommand{\const}{\rm const}
  
  \newcommand{\supp}{\rm supp}

  \newcommand{\mes}{\rm mes}

\begin{document}

   \begin{center}

 {\bf Tail estimations for functions belonging to Grand Lebesgue  }\par

\vspace{4mm}

 {\bf  Spaces builded on the set with infinite measure. } \par

\vspace{5mm}

{\bf M.R.Formica,  E.Ostrovsky and L.Sirota}.\par

\end{center}

 Universit\`{a} degli Studi di Napoli Parthenope, via Generale Parisi 13, Palazzo Pacanowsky, 80132,
Napoli, Italy. \\

e-mail: mara.formica@uniparthenope.it \\

\vspace{4mm}

Department of Mathematics and Statistics, Bar-Ilan University, \\
59200, Ramat Gan, Israel. \\

e-mail: eugostrovsky@list.ru\\

\vspace{4mm}

\ Department of Mathematics and Statistics, Bar-Ilan University, \\
59200, Ramat Gan, Israel. \\

e-mail: sirota3@bezeqint.net \\

\vspace{5mm}

\begin{center}

{\bf Abstract.}

\end{center}

\vspace{4mm}

\hspace{3mm} We establish the bilateral exact  reciprocal interrelations between a tail behavior of a measurable functions
 and its norm in the suitable Grand Lebesgue Space (GLS) as well as Orlicz one, builded over the set with {\it infinite} measure. \par
  \ We bring also some examples in order to illustrate the exactness of offered estimates. \par

\vspace{4mm}

\begin{center}

{\sc  Key words and phrases.}

\end{center}

 \ Measurability, tail of function,  ordinary  Lebesgue - Riesz  and Grand Lebesgue Spaces  (GLS)  and norms, subgaussian and anti - subgaussian
 functions,  upper and lower estimates, Young - Fenchel transform, Young - Orlicz function, Orlicz norm and space,
 examples, generating and natural generating function.\par

\vspace{5mm}

\section{Statement of problem.}

\vspace{5mm}

 \hspace{3mm} Let $ \ (X = \{x\}, \cal{M}, \mu ) \ $ be measurable  space with non - trivial sigma - finite measure $ \ \mu. \ $ \par

 \vspace{3mm}

 \ {\it  We will consider in this report only the case when } $ \ \mu(X) = \infty. \ $ \par

\vspace{3mm}

 \ Denote as usually for arbitrary measurable numerical valued function $ \ f: X \to R \ $ its Lebesgue - Riesz $ \ L(p) = L(p; X) \ $ norm

$$
||f||_p :=\left[ \ \int_X |f(x)|^p \ \mu(dx) \ \right]^{1/p}, \ 1 \le p < \infty.
$$

\vspace{3mm}

\ We recall yet another important for us definition of the tail function  $ \ T[f](t), \ t > 0  \  $ for arbitrary (measurable) function $ \ f: X \to R: \ $

\begin{equation} \label{tail def}
T[f](t) \stackrel{def}{=} \mu \{x, \ |f(x)| > t \}.
\end{equation}

\vspace{4mm}

  \ Let also $ \ \psi = \psi(p), \ p \in (a,b), \ 1 \le a < b \le \infty  \ $ be certain finite strictly positive: $ \ \inf_{p \in (a,b)} \psi(p) > 0 \ $
   numerical valued function.

\vspace{3mm}

 \   Recall the definition of the so - called Grand Lebesgue Space (GLS) norm for the function $ \ f(\cdot): \ $

\begin{equation} \label{Gpsi norm}
||f||G\psi = ||f||G\psi[a,b] \stackrel{def}{=} \sup_{p \in (a,b)} \left\{ \ \frac{||f||_p}{\psi(p)} \ \right\},
\end{equation}
and correspondent Banach  (complete) functional rearrangement invariant space $ \ G\psi = G\psi[a,b] := \{ \ f: \ ||f||G\psi < \infty \ \}. \ $ \par
 \ The set of all such a functions will be denoted by $ \ \Psi[a,b] = \{\psi(\cdot) \}; \ $ and we take the notations $ \ G\psi_0 \stackrel{def}{=} G\psi[1,\infty) \ $
 as well as

$$
G\Psi := \cup_{(a,b): \ 1 \le a < b < \infty} G\psi(a,b)
$$

and

$$
(a,b) \stackrel{def}{=} \supp  (\psi).
$$

\vspace{3mm}

 \ The function $ \ \psi = \psi(p), \ p \in (a,b) \ $ is named as ordinary {\it  generating function } for this space. \par
 \ Define formally $ \ \psi(p) = \infty \ $  for the values $ \ p \notin (a,b). \ $ \par

\vspace{3mm}

 \ The theory of these spaces is represented in many works, see e.g. \cite{anatrielloformicaricmat2016}, \cite{Buldygin-Mushtary-Ostrovsky-Pushalsky},
\cite{caponeformicagiovanonlanal2013}, \cite{Ermakov etc. 1986}, \cite{Fiorenza 1}, \cite{Fiorenza 2}, \cite{Fiorenza 3}, \cite{Fiorenza 4}, \cite{Fiorenza 5},
\cite{fioguptajainstudiamath2008}, \cite{Fiorenza-Formica-Gogatishvili-DEA2018}, \cite{fioforgogakoparakoNAtoappear}, \cite{Iwaniec},\cite{Kozachenko 1},
\cite{Liflyand}, \cite{Ostrovsky 1}, \cite{Ostr Sir CLT  mixed}, \cite{Ostrovsky 0}, \cite{Ostrovsky 2}, \cite{Ostrovsky 1}, \cite{Ostrovsky 3} etc.  As a rule
it was considered the case of finite measure: \ $ \ \mu(X) = 1. \ $  The case  $ \ \mu(X) =\infty \ $ and the operators acting in these spaces  was covered  in
particular in \cite{Kokilashvili},  \cite{Samko 1} - \cite{Samko 2}, \cite{Umarkhadzhiev 1} - \cite{Umarkhadzhiev 2}.

\vspace{5mm}

 \hspace{3mm} {\bf  We intent in this short report to  establish the bilateral exact reciprocal interrelations between tail behavior
of a measurable functions and its norm in the suitable Grand Lebesgue Space (GLS), builded over the set with infinite measure. } \par
 \ {\bf  We establish also the connections between these and Orlicz spaces.} \par

\vspace{5mm}

\hspace{2mm} {\bf Remark 1.1.} Notice that in the considered here unbounded case $ \ \mu(X) = \infty, \ $  it is very essential to
investigate the behavior of the tail function $ \ T[f](t) \ $ not only as $ \ t \to \infty, \ $ but also as $ \ t \to 0+, \ $
as long as in the general case, e.g. when the function $ \ f = f(t) \ $ is strictly positive,

$$
t \to 0+ \ \Rightarrow T[f](t) \to \mes \{\ t: |f(t)| > 0  \ \} = \mes(X) = \infty.
$$

\vspace{4mm}

 \hspace{2mm} {\bf Example 1.0.} \ The bounded case $ \ \mu(X) = 1 \ $ (probabilistic case)  was investigated in
\cite{Kozachenko 2}, \cite{Ostrovsky 1}, \cite{Ostrovsky 0}, chapters 1.2. Let for instance $ \ X = R_+ \ $ and consider the
measurable functions of the form

$$
h(x) :=  h_0(x) \ I(x \in [0,1]),
$$
where $ \ I(  x \in A)  = I(A) \ $ denotes an indicator function for the measurable set (or predicator)  $ \ A. \ $  We return to the probabilistic case. \par
 \ In particular,  suppose
\begin{equation} \label{m condit}
\exists m = \const \in (0,\infty) \ \Rightarrow  \ \forall p \in [1,\infty) \hspace{2mm}  ||h||_p \le C_1 \ p^{1/m},  \ C_1 \in (0,\infty).
\end{equation}

 \ The last relation is quite equivalent to the following tail estimate

\begin{equation} \label{m tail estim}
\exists c(m) \in (0,\infty] \ \Rightarrow  T[h](t) \le \exp \left( \  - c(m) \ t^m  \ \right), \ t \ge 0.
\end{equation}

\vspace{3mm}

 \ The case $ \ m = 2 \ $ correspondent to the famous subgaussian function (random variable.) \par

\vspace{4mm}

 \hspace{2mm} {\bf Example 1.1.} Let $ \ g: X \to R \ $ be measurable function for which $ \ \exists (a,b): \ 1 \le a < b \le \infty
 \ \forall p \in (a,b)  \Rightarrow ||g||_p < \infty.  \ $ The function

\begin{equation} \label{natural function}
\psi_g(p) \stackrel{def}{=} ||g||_p, \ p \in (a,b)
\end{equation}
is named as a { \ \it natural function} for the function $ \ g = g(x). \ $ Evidently, $ \ ||g||G\psi_g = 1. \ $\par
 \ Of course, it may be choosed as a generating function for suitable GLS; named: natural generating function. \par

\vspace{3mm}

\hspace{2mm} {\bf Sub - example 1.2.} Let $ \ X = R_+ \ $ equipped with ordinary Lebesgue measure $ \  d \mu = dx. \ $ Choose the  (positive) function

$$
g[c,\theta](x) := \exp \left(- c x^{\theta} \ \right), \ c,\theta = \const \in (0,\infty).
$$

  \ We find

$$
|| \ g[c,\theta] \ ||^p_p = \int_0^{\infty} g^p[c,\theta](x)\ dx = \int_0^{\infty} \exp \left( \ - c p x^{\theta}  \ \right) \ dx =
$$

$$
\theta^{-1} \ c^{-1/\theta} \ p^{-1/\theta} \ \Gamma(1/\theta),
$$
where as usually $ \ \Gamma(\cdot) \ $ denotes the Euler's Gamma function. So, the natural  function for $ \ g[c,\theta] \ $ has a form

$$
||g[c,\theta]||_p = \theta^{-1/p} \ (c p)^{-1/(p \theta)} \ \Gamma^{1/p} (1/\theta), \ 1 \le p < \infty.
$$

 \ A particular case: $ \ c = \theta = 1: \  g[1,1](x) = \exp(-x), \ x \ge 0; \ $

$$
||g[1,1]||_p = p^{-1/p}, \ a = 1, \ b = \infty.
$$

\vspace{3mm}

 \ Another particular case: $ \ \theta = 2, \ c = 1: \ $

\begin{equation} \label{anti subgaussian}
g[1,2](x) := \exp \left(- x^2 \right), \ x \ge 0, \ -
\end{equation}
 the so - called {\it anti - subgaussian case}.  We find

\begin{equation} \label{mom anti sub}
||g[1,2]||_p^p = 0.5 \ \pi^{1/2} \  p^{-1/2}, \ p \in (0,\infty).
\end{equation}

 \ Recall that the classical subgaussian case for the function $ \ \tilde{g} \ $ implies that

$$
||\tilde{g}||_p \asymp \sqrt{p},   \ p \ge 1.
$$

\vspace{3mm}

 \ Note that in the considered examples one can suppose $ \ p \in (0, \infty]. \ $ \par

\vspace{4mm}

 \ {\bf Remark 1.2.} It is interest in our opinion to note that in the all last examples

$$
\lim_{p \to \infty} ||g[c,\theta]||_p^p = 0,
$$
but

$$
\inf_{p \in [1,\infty)} ||g[c,\theta]||_p > 0.
$$

\vspace{4mm}

 \ {\bf Remark 1.3.} One can apply for our purpose, in the finite measure case, $ \ \mu(X) = 1 \ $ the so - called
 moment generating function

$$
\exp (\phi(\lambda)) \stackrel{def}{=} \int_X \exp(\lambda f(x)) \ \mu(dx),
$$
if this function there exists for certain non - trivial neighborhood of origin $ \ \lambda: \ |\lambda| < \lambda_0, \ $
where $ \ \exists \ \lambda_0 = \const \in (0,\infty]. \ $ \par
 \ This possibility absent in general case under considered here convention $ \ \mu(X) = \infty. \ $  In particular,

$$
  \exp (\phi(0)) = \int_X \mu(dx) = \mu(X) = \infty.
$$

\vspace{5mm}

\section{Main  result.  Upper estimate.}

\vspace{5mm}

 \ Assume that under formulated before conditions the  non - zero function $ \  f: X \to R  \ $ belongs to the certain Grand Lebesgue Space
$ \ G\psi; \hspace{2mm}  \exists \psi \in G\Psi. \ $  Introduce  the  following notations:

$$
\nu(p) = \nu[\psi](p):= p \ln \psi(p), \ p \in \supp (\psi); \ \gamma := ||f||G\psi \ \in (0,\infty);
$$

$$
\nu^*(t) = \nu^*[\psi](t) \stackrel{def}{=} \sup_{p  \in \supp (\psi)} (pt - \nu[\psi](p)) \ -
$$
the (regional) Young - Fenchel transform of the function $ \ \nu. \ $  \par

\vspace{5mm}

\ {\bf Theorem 2.1.}

\vspace{3mm}

\begin{equation} \label{upp estim}
T[f](t) \le \exp \left( \  - \nu^*(\ln (t/\gamma)) \ \right), \ t > e \ \gamma.
\end{equation}

\vspace{5mm}

 \ {\bf Proof} is quite alike to one in the finite - measure case  \cite{Kozachenko 1} - \cite{Kozachenko 2}. Indeed,
 one can suppose without loss of generality

$$
\gamma = ||f||G\psi  = 1.
$$

 \ It follows immediately from the  direct definition of the Grand Lebesgue Norm  that

$$
\forall p \in \supp (\psi) \ \Rightarrow  ||f||_p \le \psi(p),
$$
or equally

$$
\forall p \in \supp (\psi) \ \Rightarrow \int_X |f(x)|^p \ \mu(dx) \le \psi^p(p) = \exp(p \ln \psi(p)) = \exp(\nu(p)).
$$

 \ We apply the Markov - Tchebychev's inequality:

$$
T[f](t) \le \frac{\exp(\nu(p))}{t^p} = \exp \left( \ - ( p \ln t - \nu(p)) \ \right), \ t \ge e.
$$

 \ It remains to take the minimum over all the admissible values of the parameter $ \ p \in \supp (\psi). \ $  \par

\vspace{5mm}

 \ {\bf Example 2.1.} Consider the function $ \ g =  g[1,\theta](x), \ x \in R_+:  \ $

$$
g(x) = \exp  \left(-x^{\theta} \right), \theta = \const > 0.
$$

 \ The tail function for $ \ g(\cdot) \ $ has a form

$$
T[g](t) = 0, \ t > 1; \hspace{2mm} T[g](t) = |\ln t|^{1/\theta}, \ t \in (0,1).
$$
 \ The upper estimate for this function based on the theorem 2.1 and example  1.1 give
at the same up to constants result. \par

\vspace{5mm}

\section{Main result:  lower estimates.}

\vspace{5mm}

 \hspace{3mm} Given: the tail function $ \ T[f](t), \ t > 0, \ $ (or its upper estimate.)  We
 intent to  establish  under our conditions the {\it exact up to multiplicative constant}
 estimate for the GLS norm for this function  $ \ ||f||G\psi \ $ for suitable generating function
 $ \ \psi(\cdot). \ $ \par

\vspace{3mm}

 \ Note first of all that

\begin{equation} \label{mom estim}
||f||_p^p = \int_X |f(x)|^p \ \mu(dx) = p \int_0^{\infty} t^{p-1} \ T[f](t) \ dt, \ p \in (0, \infty).
\end{equation}

\vspace{3mm}

 \ Therefore, the natural generating function for the function  $ \  f(\cdot) \, \ $ namely $ \ \psi = \psi(p) = \psi[f](p) \ $
 has a form

 \vspace{4mm}

 \ {\bf Proposition 3.1.}

 \vspace{3mm}

\begin{equation} \label{nat for f}
\psi[f](p) = \left[ \   p \int_0^{\infty} t^{p-1} \ T[f](t) \ dt   \ \right]^{1/p},
\end{equation}
if of course there exists for some non - trivial segment $ \ p \in (a,b) = \supp ( \psi(p)), \ 0 < a < b \le \infty, \ $
so that

$$
f \in G\psi[f], \hspace{3mm} ||f||G\psi[f] = 1.
$$

\vspace{4mm}

 \hspace{2mm} {\bf Example 3.1.} Let  as above $ \ X = (-1,0),  \ g(x) := | \ \ln| x| \ |, \ x \in X; \ $

$$
Y := (0,\infty); \hspace{3mm}  h(y) = e^{-y}.
$$
 \ Define also

 $$
  Z = \{z\}  = \{(x,y)\} = X \otimes Y, \ X \cap Y= \emptyset;  \ f(z) =  f(x,y) := g(x) + h(y).
 $$

 \vspace{3mm}

 \ There holds

$$
||g||_p^p = \Gamma(p+1), \ ||h||^p_p  = \frac{1}{p}, \ p > 0.
$$
 \ As long as this functions  $ \ f, g \ $  are disjoint:  $ \  f(x) \cdot g(y) = 0,  \ $

$$
||f||^p_p = ||g||^p_p + ||h||^p_p
$$
 and

$$
||f||_p^p \asymp \Gamma(p+1), \ p \to \infty; \hspace{3mm} ||f||_p^p \asymp \frac{1}{p}, \ p \in (0,1).
$$

 \ The last relations stands in complete accordance with the tail behavior  of the function $ \ f: \ $

$$
T[f](t) \asymp \frac{1}{|\ln t|}, \ t \in (0,1); \ T[f](t) \asymp e^{-t}, \ t \in (1,\infty).
$$

\vspace{5mm}

\section{Orlicz space characterization of tail behavior.}

\vspace{5mm}

 \hspace{3mm} Statement of problem: given a tail function $ \ T = T(t), \ t > 0; \ $  find an Young - Orlicz function
$ \  N = N(u) = N[T](u), \ u \ge e, \ $   such that for arbitrary measurable function $ \ f: X \to R \ $ for which
$ \ T[f](t) \le T(t) \ $ this function belongs also the Orlicz space $ \ L(N) \ $ builded over $ \ (X,\cal{M}, \mu). \ $ \par
 \ We understand as the Young - Orlicz function $ \ N = N(u), \ u \in [0,\infty) \ $ the continuous non - negative strictly
 increasing function for which

$$
\lim_{u \to 0+} \frac{N(u)}{u} = 0; \hspace{3mm} \lim_{u \to \infty} \frac{N(u)}{u}  = \infty;
$$
not necessarily  to be convex. \par

 \ The finite - measure case $ \ \mu(X) = 1 \ $ is investigated in particular in \cite{Kozachenko 2}. \par

\vspace{3mm}

 \ It is convenient for us to represent the tail function $ \ T = T(t) \ $ as an exponential form

\begin{equation} \label{exp represent}
T(t) = \exp(-w(t)), \ t \ge 1.
\end{equation}

\vspace{3mm}

 \hspace{3mm} Note first of all that if for some increasing positive function $  \ G = G(t) \ $

$$
I := \int_X G(|f(x)|) \ \mu(dx) < \infty,
$$
then by virtue of Tchebychev - Markov's inequality

\begin{equation} \label{Tail estim}
T[f](t) \le I / G(t), \ t > 0 ,
\end{equation}
or equally

\begin{equation} \label{G upp estim}
G(t)  \le I /T[f](t), \ t < 0.
\end{equation}

 \ Therefore, it is reasonable to choose as the Young - Orlicz function  the tail one $ \ T = T(t): \ $

\begin{equation} \label{GT}
N[T](t) := G_0(t) = G_0[T](t) \stackrel{def}{=} \frac{1}{T(t)}, \ t \ge 1.
\end{equation}

\vspace{5mm}

 \ {\bf Theorem 4.1.}  Suppose

\vspace{3mm}

\begin{equation} \label{key condit}
\exists k = \const > 0 \ \Rightarrow \int_0^{\infty} \frac{|dT(t)|}{T(t/k)} < \infty.
\end{equation}

\vspace{3mm}

 \ Then the arbitrary measurable function $ \ f: X \to R \ $ for which $  \ T[f] (t) \le T(t), \ t > 0 \ $
 belongs to the Orlicz space $ \ L(N[T]) \ $  builded over source measurable space $ \ (X, \cal{M}, \mu). \ $  \par

\vspace{5mm}

 \ {\bf Proof.} We will use the following fact

$$
\int_X H(|f|(x)) \ \mu(dx)  = \int_0^{\infty} H(t) \ |dT[f](t)|,
$$
see e.g. \cite{Stein}, chapters 1,2; we conclude following  that for all the sufficiently greatest positive values  of
the constant $ \ K > 0 \ $

$$
\int_X N[T](|f(x)|/K) \ \mu(dx) \le \int_0^{\infty} N[T](t/K) |dT(t)| =
$$

$$
\int_0^{\infty} \frac{|dT(t)|}{T(t/K)} < \infty, \hspace{3mm}  K > k;
$$
we used a comparison inequality, see e.g. \cite{Kozachenko 2}. \par

\vspace{4mm}

 \hspace{3mm} {\bf Examples 4.1 - 4.2.} The condition (\ref{key condit}) is satisfied for the tail functions of the form

$$
\exists \ C,m \in (0,\infty) \ \Rightarrow  \  T_{m,C}(t) = \exp \left( - C t^m \ \right), \ t > 0,
$$
 and is not satisfies for the  next tail functions

$$
T^{(\theta)}(t) = |\ln t|^{1/\theta} \cdot I(t \in (0,1)), \ \theta \in (0,\infty).
$$

 \ The finding of the correspondent Orlicz function for this tail function  is an open problem.\par

\vspace{5mm}

\section{Concluding remarks.}

\vspace{5mm}

 \hspace{3mm}  {\bf  \ Possible generalization.}  It  is interest, in our opinion, to generalize  the obtained results on the
 multidimensional case, i.e. when the functions $ \ \{f\} \ $  are vector valued, in the spirit of the preprint \  \cite{Ostrovsky Mult}. \par

\vspace{6mm}

\vspace{0.5cm} \emph{Acknowledgement.} {\footnotesize The first
author has been partially supported by the Gruppo Nazionale per
l'Analisi Matematica, la Probabilit\`a e le loro Applicazioni
(GNAMPA) of the Istituto Nazionale di Alta Matematica (INdAM) and by
Universit\`a degli Studi di Napoli Parthenope through the project
\lq\lq sostegno alla Ricerca individuale\rq\rq .\par


\begin{thebibliography}{44}


\bibitem{anatrielloformicaricmat2016}
{\bf G.~Anatriello} and {\bf M.~R.~Formica.} {\it Weighted fully
measurable grand Lebesgue spaces and the maximal theorem}. Ric. Mat.
\textbf{65} (2016), no.~1, 221--233.


\bibitem{Buldygin-Mushtary-Ostrovsky-Pushalsky}
{\bf V.~V.~Buldygin, D.~I.~Mushtary, E.~I.~Ostrovsky} and {\bf M.~I.~Pushalsky.} {\it New Trends in
Probability Theory and Statistics.} Mokslas (1992), V.1, 78--92;
Amsterdam, Utrecht, New York, Tokyo.



\bibitem{caponeformicagiovanonlanal2013}
{\bf C.~Capone, M.~R.~Formica} and {\bf R.~Giova.} {\it Grand
{L}ebesgue spaces with respect to measurable functions}. Nonlinear
Anal. \textbf{85} (2013), 125--131.


\bibitem{Ermakov etc. 1986}
{\bf S. V. Ermakov, and E. I. Ostrovsky.} {\it Continuity Conditions, Exponential Estimates, and the Central Limit Theorem for Random Fields.}
 Moscow, VINITY,  1986. (in Russian).



\bibitem{Fiorenza 1}
{\bf A. Fiorenza, G. E. Karadzhov.}  {\it Grand and small Lebesgue spaces and their analogs.}
Journal for Analysis and its Applications, 23(4) (2004), 657 \ - \ 681.


\bibitem{Fiorenza 2}
{\bf A. Fiorenza, J. M. Rakotoson.}  {\it Petits espaces de Lebesgue et leurs applications.}
C.R.A.S. t, 333 (2001), 1–4.

\bibitem{Fiorenza 3}
{\bf A. Fiorenza C. Capone and G.E. Karadzhov.} {\it  Grand Orlicz spaces and global integrability
of the Jacobian.}  Math. Scand., {\bf 102(1),} 131 \ - \ 148, 2008.


\bibitem{Fiorenza 4}
{\bf A. Fiorenza.}  {\it Duality and reflexivity in grand Lebesgue spaces.} Collect. Math.,
{\bf 51(2)}, 131 \ - \ 148, 2000.

\bibitem{Fiorenza 5}
{\bf C. Capone and A. Fiorenza.}  {\it On small Lebesgue spaces.} J. Function Spaces and
Applications, {\bf 3,} 73 \ - \ 89, 2005.



\bibitem{fioguptajainstudiamath2008}
{\bf A.~Fiorenza, B.~Gupta} and {\bf P.~Jain.} {\it The maximal
theorem for weighted grand Lebesgue spaces}. Studia Math.
\textbf{188} (2008), no.~2, 123--133.



\bibitem{Fiorenza-Formica-Gogatishvili-DEA2018}
{\bf A.~Fiorenza, M.~R.~Formica} and {\bf A.~Gogatishvili.} {\it On
grand and small Lebesgue and Sobolev spaces and some applications to
PDE's}. \emph{Differ. Equ. Appl.} \textbf{10} (2018), no.~1, 21--46.



\bibitem{fioforgogakoparakoNAtoappear}
{\bf A.~Fiorenza, M. R.~Formica, A.~Gogatishvili, T.~Kopaliani} and
{\bf J.~M. Rakotoson.} {\it Characterization of interpolation
between grand, small or classical Lebesgue spaces}. Preprint
arXiv:1709.05892, Nonlinear Anal.




\bibitem{Iwaniec}
{\bf T.Iwaniec, C. Sbordone.}  {\it On the integrability of the Jacobian under minimal hypotheses. }
Arch. Rational Mech. Anal., {\bf 119,}  (1992), 129 \ - \ 143.


\bibitem{Kokilashvili}
{\bf V.Kokilashvili.} {\it  Weighted problems for operators of harmonic analysis in some
Banach function spaces.}  Lecture course of Summer School and Workshop
"Harmonic Analysis and Related Topics" (HART2010), Lisbon, June 21-25,
http://www.math.ist.utl.pt/~hart2010/kokilashvili.pdf, 2010.

\bibitem{Kozachenko 1}
{\bf Kozachenko Yu. V., Ostrovsky E.I.} (1985). {\it The Banach Spaces of random
Variables of subgaussian type.} Theory of Probab. and Math. Stat. (in
Russian). Kiev, KSU, 32, 43 - 57.


\bibitem{Kozachenko 2}
{\bf Kozachenko Yu.V., Ostrovsky E., Sirota L.}
{\it Equivalence between tails, Grand Lebesgue Spaces and Orlicz
norms for random variables without Cramer's condition.} \\
arXiv:1710.05260v1 [math.PR] 15 Oct 2017

\bibitem{Liflyand}
{\bf E. Liflyand, E. Ostrovsky, L. Sirota.} {\it Structural properties of Bilateral Grand
Lebesque Spaces.} Turk. J. Math., {\bf 34,}  (2010), 207 \ - \ 219.

\bibitem{Ostrovsky 1}
{\bf Ostrovsky E. and Sirota L.}  {\it Moment Banach spaces: theory and applications.}
HIAT Journal of Science and Engineering, C, Volume {\bf 4,}  Issues 1-2, pp. 233 \ - \ 262, (2007).

\bibitem{Ostr Sir CLT  mixed}
{\bf E.Ostrovsky, L.Sirota.}   {\it Central Limit Theorem and exponential estimations in mixed (anisotropic)
Lebesgue Spaces.} \par
arXiv:1308.5606v1 [math.PR] 26 Aug 2013

\bibitem{Ostrovsky 0}
{\bf Ostrovsky E.I.}  (1999). {\it Exponential estimations for random Fields and its
Applications, (in Russian).}  Moscow - Obninsk, OINPE.

\bibitem{Ostrovsky 2}
{\bf E.Ostrovsky, L.Sirota and E.Rogover.} {\it Integral operators in bilateral Grand Lebesgue Spaces.} \par
arXiv:0912.2538v1 [math.FA] 13 Dec 2009

\bibitem{Ostrovsky Mult}
{\bf Ostrovsky E., Sirota L.}  {\it Some new rearrangement invariant spaces: theory and
applications.} Electronic publications: arXiv:math.FA/0605732 v1, 29, (May 2006);

\bibitem{Ostrovsky 3}
{\bf E.Ostrovsky, L.Sirota} {\it Central Limit Theorem and exponential tail estimations in mixed
(anisotropic) Lebesgue spaces.} \\
arXiv:1308.5606v1 [math.PR] 26 Aug 2013

\bibitem{Samko 1}
{\bf S.G. Samko, S.M. Umarkhadzhiev.}  {\it On Iwaniec - Sbordone spaces on sets which may
have infinite measure.}  Azerb. J. Math. {\bf 1, } (1), (2010),  67 \ - \ 84.


\bibitem{Samko 2}
{\bf S.G. Samko, S.M. Umarkhadzhiev.}  {\it On Iwaniec - Sbordone spaces on sets which may have infinite measure.}
Azerbaijan Journal of Mathematics, V. {\bf 1,}  No 1, 2011, January, ISSN 2218 \ - \ 6816.

\bibitem{Stein}
{\bf E.M.Stein.}  {\it Singular Integrals and Differentiability Properties of Functions.}
Princeton Univ. Press, 1970, Princeton, New York.

\bibitem{Umarkhadzhiev 1}
{\bf Umarkhadzhiev, S. M.}  {\it Generalization of the Notion of Grand Lebesgue Space.}  Russian Mathematics,
2014, vol. 58, no. 4, pp. 35 \ - \ 43. DOI: 10.3103/S1066369X14040057.



\bibitem{Umarkhadzhiev 2}
{\bf Umarkhadzhiev, S. M.}  {\it  Integral Operators with Homogeneous Kernels in Grand Lebesgue Spaces.}
Mathematical Notes, 2017, vol. 102, no. 5 \ - \ 6, pp. 710 \ - \ 721. DOI: 10.1134/S0001434617110104.




\end{thebibliography}
\end{document}